\documentclass[a4paper]{article}

\usepackage{amsmath}
\usepackage{amssymb}
\usepackage[all]{xy}

\newtheorem{theo}{Theorem}[section]

\newtheorem{prop}[theo]{Proposition}

\newtheorem{lem}[theo]{Lemma}
\newtheorem{rem}[theo]{Remark}
\newtheorem{exem}[theo]{Example}

\newtheorem{conj}{Conjecture}
\newtheorem{theorem}[conj]{Theorem}
\newtheorem{corol}[conj]{Corollary}
\newtheorem{proposition}[conj]{Proposition}

\renewcommand{\P}{\mathbf{P}}
\newcommand{\C}{\mathbf{C}}

\newcommand{\Q}{\mathbf{Q}}
\newcommand{\Z}{\mathbf{Z}}
\renewcommand{\O}{\mathcal{O}}
\renewcommand{\H}{\mathrm{H}}

\newcommand{\Spec}{\mathrm{Spec}}
\newcommand{\Pic}{\mathrm{Pic}}
\newcommand{\phitilde}{\widetilde{\varphi}}
\newcommand{\Stilde}{\widetilde{S}}
\renewcommand{\epsilon}{\varepsilon}
\renewcommand{\hat}{\widehat}
\newcommand{\M}{\mathcal{M}}
\newcommand{\V}{\mathcal{V}}

\newcommand{\para}{\refstepcounter{theo} 
\vspace{0.2cm} \noindent
\textbf{\thetheo~}}

\newcounter{paragraf}[section]

\author{Thomas Dedieu\\
Institut de math{\'e}matiques de Jussieu, UMR 7586}
\title{Severi varieties and self rational maps of $K3$ surfaces}

\begin{document}

\maketitle

\section*{Introduction}

\para\label{notations} \textbf{Notations.} 
We deal in this paper with complex projective $K3$ surfaces, \emph{i.e.}
smooth $K$-trivial complex projective surfaces without irregularity. 
Let $\varphi : S \dashrightarrow S$ be a dominant self
rational map.  Suppose $\Pic(S)=\Z$. Then there exists a positive
integer $l$ such that
$\varphi^* \O_S(1) \cong \O_S(l)$. It is the 
algebraic degree of $\varphi$, that is the degree of the polynomials
defining $\varphi$. There always exists an elimination of
indeterminacies 
$$\xymatrix{
\widetilde{S} \ar[d]_{\tau} \ar[dr]^{\widetilde{\varphi}} \\
S \ar@{-->}[r]_{\varphi} & S
}$$
\emph{i.e.} a commutative diagram, where $\phitilde$ is a morphism and
$\tau$ is a finite sequence of blow-ups. One defines the topological degree of
$\varphi$ as $\deg \varphi := \deg \phitilde$. It is the number of points in
the inverse image of a generic point $x \in S$ under the action of
$\varphi$. We write $R$ for the
ramification divisor of $\phitilde$ (it is the zero divisor of the
Jacobian $\bigwedge^2 d\phitilde$).

\para \textbf{Self-rational maps.}
The main goal of this article is to study the geometric and
numerical properties of self-rational maps $S \dashrightarrow S$,
where $S$ is a $K3$ surface with Picard group $\Z$, in order to attack
the following conjecture.
\begin{conj}\label{endom}
For a generic projective $K3$ surface $S$, there does not exist any
dominant rational map
$$\varphi : S \dashrightarrow S$$
satisfying $\varphi^* \O_S(1) \cong \O_S(l)$, $l>1$.
\end{conj}

It is another problem to find complex projective manifolds $X$
equipped with a self morphism  $f : X \rightarrow X$. It has already
been studied by Beauville (\cite{ba}), Fujimoto and Nakayama
(\cite{fujimoto}, \cite{nakayama}, \cite{fuji-naka}), or Amerik,
Rovinsky and Van de Ven (\cite{arv}). 
Toric and abelian varieties are obvious examples of such manifolds,
and one does not know any other example
than those deduced from these two obvious ones. This leads us to
conjecture that there does not exist any other example at all, or in
other words that in case $\Pic(X)=\Z$ and $\kappa(X)\geqslant 0$,
there does not exist any $f$ with $\deg f >1$.
Beauville (\cite{ba}) proves in this direction that a complete
intersection of $p$ hypersurfaces of respective degrees
$d_1,\ldots,d_p$ in $\P^{n+p}$ ($n \geqslant 2$) does not admit any
endomorphism of degree strictly greater than 1, as soon as at least
one of the $d_i$'s is greater than 3. 
One easily sees that a $K3$ surface cannot possess any dominant
endomorphism of degree strictly larger than 1. Indeed, 
if $S$ is a $K3$ surface and $f : S \rightarrow S$ a dominant
morphism, then the relation
$$K_S = f^* K_S + R$$
(where $R$ is the ramification divisor of $f$), combined with the fact
that $K_S=0$, proves that $f$ is necessarily smooth. It is thus an
{\'e}tale cover of $S$ by itself, and therefore an automorphism of $S$,
because all $K3$ surfaces are simply connected.

Amerik and Campana (\cite{ac}), or Cantat (\cite{cantat}) are more
generally interested in the search of manifolds equipped with dominant
self rational maps.
Again, we have easy examples deduced from the toric and abelian
ones. If $S$ is a $K3$ surface equipped 
with an elliptic pencil $|F|$ and a relatively ample line bundle $L$
of relative degree $d$,
then we can construct $\mu_{d+1} : S \dashrightarrow S$ of degree
$(d+1)^2$, \emph{via} multiplication in the fibers of
the pencil~: the image of a point $x$ on a generic fiber $F$ is
defined as the unique point $y \in F$ satisfying 
$$\O_F\left( (d+1)x - y \right) = \left. L \right|_F.$$
Kummer surfaces are smooth models of quotients of complex tori under
the action of an involution. One can therefore construct self
rational maps of degree strictly greater than 1 on them, by descending
the homotheties on the tori. Note that these two examples only concern
special $K3$ surfaces.

In greater dimensions, one has examples that cannot be deduced from
the two obvious ones. Voisin
constructs in \cite{Kcorresp} a self rational map of degree 16 of
the variety $X$ of lines in a cubic hypersurface $V \subset \P^5$. 
$X$ is 4-dimensional 
and hyperk{\"a}hler, has $\Pic(X) \cong \Z$, and is deformation equivalent
to the punctual Hilbert scheme $S^{[2]}$ of a $K3$ surface $S$
of degree 14 (\cite{BeauDo}). One gets this self rational map by
mapping the generic line 
$l \subset V$ to the residual line $l'$ to $l$ in $P\cap V$, where $P$ is
the unique 2-plane tangent to $V$ along $l$. This map does not respect
any fibration in virtue of the following theorem of Amerik and Campana
(\cite{ac})~: 
if $X$ is a projective manifold satisfying $K_X=0$ and
  $\mathrm{NS}(X)=\Z$, then any rational fibration $g : X
  \dashrightarrow B$ ($0< \dim B < \dim X$) has fibers of general type.
Voisin's example shows in particular that conjecture \ref{endom} is
specific to the case of surfaces.

The holomorphic dynamical point of view gives a new insight into this
problem. Given a complex
manifold $X$ equipped with a transformation $f : X \rightarrow X$, one
gets a discrete dynamical system by iterating $f$ (see
\cite{cantat}). Amerik and Campana 
associate a meromorphic fibration $g : X \dashrightarrow T$ to any
dominant self rational map $f : X \dashrightarrow X$. Its general
fiber $X_t$ is the Zariski closure of the orbit of a general point in
$X_t$. This allows them to prove that if $X$ is a complex projective
manifold satisfying both $K_X=0$ and $\mathrm{NS}(X) \cong \Z$, and
if $f : X \dashrightarrow X $ is a dominant self rational map of
degree greater than 2, then the orbit of a general point of $X$ under
the action of $f$ is Zariski-dense. 

Such dynamical facts induce results concerning potential density in
the $K$-trivial case (a
variety $X$ over some field $k$ is said to be potentially dense if
there exists some finite extension $k \rightarrow k'$, such that the
set of $k'$-rational points is Zariski-dense in $X$, see \emph{e.g.}
\cite{ht}). 
As a simple corollary of their results exposed above, Amerik and
Campana get : if $X$ is a smooth projective variety defined over
$\bar{k}$ ($k$ a non countable field), such that $K_X=0$ and
$\Pic(X)=\Z$, and if there exists $f : X \dashrightarrow X$ with $\deg
f >1$, then $X$ has potential density.
Cantat gets on his side, and using dynamical methods, a large part of
the proof of the following theorem, due to Bogomolov and Tschinkel 
(\cite{bt})~: a projective $K3$ surface defined over some number field
$k$ is potentially dense, as soon as it can be realized as an
elliptic fibration. 

Eventually, this dynamical study allows Cantat to show that if $X$ is
a projective Calabi-Yau manifold of dimension $n$, and if there exists
$f : X \dashrightarrow X$ satisfying some dilating property, then
there exists a dominant rational map $\C^n \dashrightarrow X$ (in
particular, $X$ satisfies the Kobayashi conjecture, see \emph{e.g.}
\cite{harvard}).  
This leads him to ask the following questions concerning a generic
algebraic $K3$ surfaces $S$ (in addition to the question 
of the existence of a dominant self-rational map $\varphi : S
\dashrightarrow S$ with $\deg \varphi > 1$, to which this article
gives a conjectural answer)~: does $S$ have potential density~? does
$S$ admit a dominant rational map $\C^2 \dashrightarrow S$~?

\para \textbf{Severi varieties.} 
We present in this article a result relating conjecture~\ref{endom}
and the irreducibility of Severi varieties for $K3$ surfaces.

Nodal plane curves (\emph{i.e.} plane curves with only non-degenerate
singularities) are a classical topic. A historical reason for this is
the fact that every smooth curve is birationally equivalent to a nodal
plane curve, \emph{via} a series of projections. Let $V_{d,g}$ be the
variety parametrizing plane irreducible curves of degree $d$ and
geometric genus $g$. It is called a Severi variety. It is the closure
in the projective space parametrizing all plane curves of degree $d$
of the locus of irreducible nodal curves of genus $g$. Severi gave an
uncomplete proof of the fact that all varieties $V_{d,g}$ are
irreducible, and it is only in 1986 that Harris actually proved this
(see \cite{severi}).

A natural generalization of this is the study of nodal curves on a
projective surface $S$ equipped with a fixed ample effective line
bundle $L$. $V_{k,h}$ then denotes the closure in $|kL|$ of the locus
of irreducible nodal curves of geometric genus $h$, and is again called
a Severi variety. The following questions arise naturally~: when are
the $V_{k,h}$ non empty ? What are their dimensions~? Are they
irreducible~? smooth~?
They are studied by Chiantini and Ciliberto in \cite{ciliberto}. Also
Greuel, Lossen and Shustin (\cite{gls}), Keilen (\cite{keilen}) give
some numerical criteria for generalized Severi varieties to be
irreducible.

We focus here on universal Severi varieties for $K3$ surfaces. Fix an
integer $g \geqslant 2$, and write $\M_{K3,g}$ for the moduli space of
$K3$ surfaces equipped with an indivisible, ample line bundle $L$ of
self-intersection $2g-2$ (we call these $K3$ surfaces of genus
$g$). There exists a universal family $\mathcal{S}_g \rightarrow
\M^{\circ}_{K3,g}$ over an open subset of $\M_{K3,g}$. To a generic
point $m \in \M^{\circ}_{K3,g}$ corresponds a $K3$ surface $S_m$
with Picard group
$$ \Pic(S_m) = \Z \cdot L_m, $$
where $L_m$ is an ample and indivisible divisor class, satisfying
$L_m^2 = 2g-2$ (see \emph{e.g.} \cite{palaiseau}). 
A generic member of the complete
linear system $|L_m|$ is a smooth curve of geometric genus $g$. For
integers $k,h \geqslant 1$, we define the universal Severi variety 
$$ \V_{k,h} \longrightarrow \M^{\circ}_{K3,g} $$ 
to be the variety whose fiber over a generic $m \in \M^{\circ}_{K3,g}$
is the closure in $|kL_m|$ of the locus 
$$ \left\{
C \in |kL_m|\ \mbox{s.t.}\ C\ \mbox{is irreducible, nodal, and of
  geometric genus } h
\right\}. $$
By the genus formula, all curves in the complete linear system
$|kL_m|$ have arithmetic genus $p_a(k)=1 + (kL_m)^2/2 = 1+k^2(g-1).$

The deformation theory of nodal curves on $K3$ surfaces works very
well. In particular, we know that for a $K3$ surface $S$, the Severi variety
$V_{k,h}$ is smooth and of the expected dimension~; if $S$ is generic,
then $V_{k,h}$ is non empty (see section
\ref{familles}). The only question that remains open is whether the
$V_{k,h}$ are irreducible or not. For $h=0$, it is clear that the
answer is no~: there are finitely many rational curves in the linear
system $|kL|$. The corresponding Severi variety is then a disjoint
union of points, which of course is not irreducible. 

It is perfectly possible that the universal Severi variety $\V_{k,h}$ is
irreducible even if the fibers $V_{k,h}$ are reducible. The question
of the irreducibility of the $\V_{k,h}$ is the closest to the initial
problem of Severi, where the projective plane plays the role of a
universal space for complete non singular curves, since they all are 
birationally equivalent to plane curves with at most nodes as
singularities. 
We conjecture that all $\V_{k,h}$ are irreducible. The following
less optimistic version is however sufficient for our purpose.
\begin{conj}\label{severi}
Let $\epsilon >0$ be given. If $k$ is great enough with regard to
$\epsilon$, then for all integer $h$ satisfying
$$ \epsilon p_a(k) \leqslant h \leqslant p_a(k), $$
the universal Severi variety $\V_{k,h}$ is irreducible.
\end{conj}

\para \textbf{Results.}
In section \ref{lien}, we prove the following.
\begin{theorem}
Let $g,l \geqslant 2$ be given. If for $m \in \M^{\circ}_{K3,g}$
generic there exists
a dominant rational map $\varphi_m : S_m \dashrightarrow S_m$
satisfying $\varphi_m^* \O_S(1) \cong \O_S(l)$, 
then for $k$ great enough the universal Severi variety
$\mathcal{V}_{kl,p_a(k)}$ possesses at least two irreducible
components. 
\end{theorem}

To prove this, we look at the images under the action
of $\varphi$
of the curves in $|kL|$. We show 
that they are elements of $|klL|$, and are generically nodal and of
geometric genus $p_a(k)$. This gives a way to construct two
distinct irreducible components of $\V_{kl,p_a(k)}$, the first one
parametrizing curves whose respective singularity $0$-cycles are all
rationally equivalent to a constant, and the second one parametrizing
curves with non constant singularity $0$-cycle modulo rational
equivalence. By singularity
$0$-cycle of a curve $C \subset S$, we mean the sum of all singular
points of $C$, seen as a $0$-cycle on $S$.

Now asymptotically, we have $p_a(k)/p_a(lk) \sim_{k
\rightarrow \infty} 1/l^2.$  We thus get the following result, which
gives a way to tackle conjecture~\ref{endom}.
\begin{corol}
Conjecture \ref{severi} on the Severi varieties implies conjecture
\ref{endom} on self-rational maps.
\end{corol}

In section \ref{contraintes}, we gather numerical constraints between
the topological and numerical degrees (\emph{i.e.} $\deg \varphi$ and
$l$ with the notations of \ref{notations}) of a dominant self-rational map
$S \dashrightarrow S$, where $S$ is a given generic $K3$ surface. This
restricts the possibilities for the existence of such self-rational
maps, and may lead to special cases of conjecture \ref{endom}. The
most significant results we get in this direction are the following.
\begin{theorem} Let $S$ be a $K3$ surface of genus $g$, with $\Pic(S)
  = \Z$. We assume there exists a dominant self rational map $\varphi
  : S \dashrightarrow S$ with $\deg \varphi >1$. Then we have the
  following. \\
(i) There exists an integer $\lambda$, such that
$ \deg \varphi = \lambda^2.$ In addition, $2g-2$ necessarily divides
$l-\lambda$ (note that we do not know the sign of $\lambda$). \\ 
(ii) There exist positive integers $\beta_1,\ldots,\beta_p$, such that
$$ l^2 = \deg \varphi + (2g-2) \textstyle{\sum_{i}} \beta_i^2. $$
$\sum_i \beta_i$ is always divisible by $2$. 
If we can eliminate the indeterminacies without any chain of
successive blow-ups of length strictly larger than $2$,
then 
$$ \deg \varphi \leqslant 1 + \frac{1}{24}\left[ p + 4(g-1) 
\textstyle{\sum_i}  \beta_i \right]. $$
\end{theorem}

We prove this by studying the geometry of an elimination of
indeterminacies of $\varphi$. To do this, we use the intersection tree
of the irreducible exceptional curves that appear in the elimination
of indeterminacies. 
We get on our way the following result,
which allows us to control the complexity of such an elimination of
indeterminacies. 
\begin{proposition}
(i) The depth of the tree is less
than $\deg \varphi -2$. If this maximal depth is achieved, then all
exceptional curves project on the same point of $S$. \\
(ii) If the tree has two connected components of
depths $l_1$ and $l_2$, then $l_1 + l_2 \leqslant \deg \varphi - 2$.\\
(iii) If all irreducible exceptional curves are disjoint, or
equivalently, if the indeterminacy can be solved with one blow-up, 
then there are at most $8(\deg \varphi -1)$ such irreducible
exceptional curves.
\end{proposition}

\para \textbf{Acknowledgements.} 
I wish to thank Claire Voisin for introducing me to this subject and for
sharing her ideas with me. Only her constant help made it possible for
me to complete this work.

\section{Families of nodal curves on $K3$ surfaces}\label{familles}

We recall here classical results on families of nodal curves in
general, and on families of nodal curves on $K3$ surfaces in
particular. We use a description with infinitesimal deformations, and
refer to \cite{voisin} (chapter 14) for more details (see
\cite{tannenbaum} for a slightly different point of view).

Let $X$ be a smooth, $n$-dimensional variety, equipped with an ample
line bundle $L$. We write $X_{(\delta)}$ for the open subset of the
symmetric product $\mathrm{Sym}^{\delta} X$ corresponding to sums $x_1
+ \cdots + x_{\delta}$ of $\delta$ pairwise distinct points. We look
at hypersurfaces in the complete linear system $|L|$ with at least
$\delta$ singular points, and therefore introduce the incidence variety
$$ \mathcal{I} = \{ (D,x_1+\cdots+x_{\delta}) \in |L| \times
X_{(\delta)} \mbox{ s.t. }D\mbox{ is singular at
}x_1,\ldots,x_{\delta} 
\}.$$ 
Let $\pi : \mathcal{I} \rightarrow |L|$ be the natural projection. In
the neighbourhood of a point $(f,x) \in  |L| \times X_{(\delta)}$
(where $f$ is an equation for a hypersurface $D$ and $x=x_1 + \cdots +
x_{\delta}$), $\mathcal{I}$ is defined by the $\delta (n+1)$ equations 
$$ f(x_i) = \frac{\partial f}{\partial z_1}(x_i) = \cdots = 
\frac{\partial f}{\partial z_n}(x_i) = 0$$
$(1 \leqslant i \leqslant \delta), $
where $z_1,\ldots,z_n$ denote with a little abuse of notations local
holomorphic coordinates at the neighbourhood of every $x_i$. We thus
have 
$$ \dim(\mathcal{I}) \geqslant \dim |L| - \delta. $$
$\dim |L| - \delta$ is called the expected dimension of
$\mathcal{I}$. 

After differentiation, we get the equations of the tangent space
$T_{\mathcal{I},(f,x)}$ in $T_{|L|,f} \times T_{X_{(\delta)},x}$~:
$(g,h_1+\cdots+h_{\delta}) \in T_{|L|,f} \times T_{X_{(\delta)},x}$
lies in $T_{\mathcal{I},(f,x)}$ if and only if one has for 
every $x_i$
$$\left\{ \begin{array}{l}
df_{x_i}(h_i) + g(x_i) = 0 \\
d \left( \frac{\partial f}{\partial z_j} \right)_{x_i} (h_i) +
\frac{\partial g}{\partial z_j} (x_i) = 0 \quad (1 \leqslant j
\leqslant n),
\end{array} \right.$$
that is
$$\left\{ \begin{array}{l}
g(x_i) = 0 \\
\mathrm{Hess}_{x_i}(f)(h_i)
= -\left(\frac{\partial g}{\partial z_1} (x_i),\ldots,\frac{\partial
  g}{\partial z_n} (x_i)\right),
\end{array} \right.$$
since all differentials $df_{x_i}$ vanish. 
The kernel of the differential $\pi_*$ at the point $(f,x)$ is 
given by $g=0$, so $\pi$ is an imbedding
at the point $(x_1+\cdots+x_{\delta},D)$ if and only if 
$\mathrm{Hess}_{x_1}(f),\ldots,\mathrm{Hess}_{x_{\delta}}(f)$ are non
degenerate, \emph{i.e.} if and only if the $\delta$
points $x_1,\ldots,x_{\delta}$ are non degenerate singular points of
$D$. In this case, the image
of $\pi_*$ at the point $(f,x)$ is simply
$$ \{ g \in T_{|L|,f} \mbox{ s.t. } g(x_1)= \ldots = g(x_{\delta}) = 0
\},$$
and the tangent space to the projection of $\mathcal{I}$ on $|L|$ at
$D$ identifies with
$$\H^0(X,\O_X(D) \otimes I_x) / \H^0(X,\O_X),$$
where $I_x \subset \O_X$ is the ideal sheaf defining
$x$. $\mathcal{I}$ is of the expected dimension $\dim |L| - \delta$ if
and only if the $\delta$ non degenerate singular points impose
independent conditions on the linear system $|L|$. In this case,
the non degenerate singular points of $D$ can be independently smoothed
by deformation.

\begin{theo}\label{dimfamille}
Let $S$ be a generic $K3$ surface, and $L$ an ample and indivisible
line bundle on it. For any positive integer $k$, and any $h \leqslant
p_a(k)$, the quasi-projective variety
$$ V_{k,h}^{\circ} = \{ C \in |k L| \mbox{ s.t. } C \mbox{ is
  irreducible, nodal, and of geometric genus }h \} $$
is non empty, smooth, and of the expected dimension $h$.
\end{theo}

\noindent \textbf{Proof.}
Riemann-Roch formula for surfaces gives $\dim |kL| = p_a(k)$, so $h$ is
the expected dimension of $V_{k,h}$. Suppose we have some irreducible
nodal curve $C \in |kL|$ with precisely $\delta=p_a(k) - h$ nodes. Let
$Z$ denote the singularity $0$-cycle of $C$. The former
infinitesimal calculations give 
the Zariski-tangent space
$$ T_{V_{k,h},C} \cong \H^0(S,\O_S(C) \otimes I_Z) /
\H^0(S,\O_S) \cong \H^0(C,\O_C(C) \otimes I_Z).$$ 
The canonical bundle $K_S$ being trivial, we get by adjunction formula
$$ \O_C(C) \cong \O_C(K_S+C) \cong \omega_C.$$
Since $C$ is nodal, we have $K_{\widetilde{C}} = \nu^* K_C(-2Z)$,
where $\nu : \widetilde{C} \rightarrow C$ is a 
normalization of $C$. Therefore
$$ \H^0(C,\O_C(C) \otimes I_Z) \cong \H^0(C,\nu_* 
K_{\widetilde{C}}), $$ 
which gives $\dim T_{V_{k,h},C} = g(C) = h$. Since $\dim V_{k,h}
\geqslant \dim |kL| - \delta = h$, this proves that $V_{k,h}^{\circ}$ is
smooth and of the expected dimension.

It therefore only remains to show that there actually exists an
irreducible nodal curve $C \in |kL|$ with precisely $\delta=p_a(k) -
h$ nodes. It is enough to find an irreducible rational nodal curve in
$|kL|$, since it gives a genus $h$ curve by smoothing exactly $h$ of
its nodes. This is given by Chen's theorem below.

\hfill $\Box$

\begin{theo}[Chen, \cite{chen}]
Consider $n\geqslant 3$ and $k>0$. For $S$ a generic $K3$ surface in
$\P^n$, the complete linear system $|\O_S(k)|$ contains an irreducible
rational curve with only nodes as singularities.
\end{theo}

\section{Link between conjectures \ref{endom} and \ref{severi}}\label{lien}

In this section, we prove that conjecture~\ref{severi} implies
conjecture~\ref{endom}. We start with the following result, which
for a generic $K3$ surface $S$, describes the geometric action of a
dominant self-rational map $\varphi : S \dashrightarrow S$ on a
generic curve $C \in |\O_S(k)|$.

\begin{prop}\label{action}
Let $S$ be a $K3$ surface of genus $g \geqslant 2$, with $\Pic(S)=\Z$,
and assume there exists a dominant rational map $\varphi
: S \dashrightarrow S$ satisfying $\varphi^*\O_S(1) \cong \O_S(l)$. We
consider $C \in |\O_S(k)|$ generic. \\
(i) Its image $\varphi(C)$ lies in $|\O_S(kl)|$. \\
(ii) For $k$ big enough, $\varphi(C)$ is irreducible and nodal, and
$C$ and $\varphi(C)$ have the same geometric genus $p_a(k)$.  
\end{prop}

\noindent \textbf{Proof.}
\emph{(i)} Consider 
$$\xymatrix{
\widetilde{S} \ar[d]_{\tau} \ar[dr]^{\widetilde{\varphi}} \\
S \ar@{-->}[r]_{\varphi} & S
}$$
an elimination of indeterminacies of $\varphi$, and write $R$ for the
ramification divisor of $\phitilde$. Since $K_{\widetilde{S}} =
\widetilde{\varphi}^* K_S + R$ and $K_S$ is trivial, $R$ is entirely
exceptional. In other words, $\varphi$ is smooth away from the
indeterminacy locus. 

Since $C \in |\O_S(k)|$ is generic, we can assume
that it avoids the indeterminacy locus. Then $\varphi_{|C}$ is locally
an imbedding, and in particular we have the equality of homology classes
$$ \left[ \varphi(C) \right] = \varphi_* [C] $$
($\varphi_*$ and $\varphi^*$ are defined as $\phitilde_* \tau^*$ and
$\tau_* \phitilde^*$ respectively).
We then compute the intersection product
$$\left< \left[ \varphi(C) \right], L \right>
= \left< \varphi_* [C],L \right> 
= \left< [C], \varphi^* L \right>,$$
where $L$ is the divisor class corresponding to $\O_S(1)$. Finally, we
have $\left< \left[ \varphi(C) \right], L \right> = kl(2g-2)$, and
therefore $\varphi(C) \in |klL|$.

\emph{(ii)} We define a scheme $S \times_S S$, which is pointwise the
set of pairs 
of points in $S$ having the same image under $\varphi$, by considering
a morphism $\varphi_U : U \rightarrow S$ ($U$ Zariski-open subset of
$S$) representing the rational map $\varphi$, and taking $S \times_S S$
to be the Zariski-closure of $U \times_S U$ in $S \times_{\Spec \C}
S$.

We claim that for $k$ big enough, $\varphi_{|_C}$ is everywhere injective
but at a finite number of points of $C$, or equivalently that 
$$C\times_S C \subset S \times_S S$$
only possesses a finite number of points outside from the diagonal. To
prove this, we define the incidence variety
$$ J = \{ (C,x_1+x_2) \in |\O_S(k)| \times S_{(2)} \mbox{ s.t. } (x_1,x_2)
\in C \times_S C \},$$
which parametrizes the pairs of distinct points of $S$ having the same
image by $\varphi$. It is given by the equations
$$\left\{ \begin{array}{l}
x_1,x_2 \in C \\
x_1 \not= x_2 \\
\varphi(x_1) = \varphi(x_2).
\end{array} \right. $$
Now the projection of $J$ on $S_{(2)}$ is $S\times_S S$, which is
pointwise the set of sums $x_1+x_2$ with 
$x_1\not=x_2$ and $\varphi(x_1) = \varphi(x_2)$, and is of dimension
2. When $k$ is large enough, the
fibers of $J$ over its projection on $S_{(2)}$ are of 
dimension $\dim |\O_S(k)|-2$ (see \ref{condind}), so 
$$ \dim J = \dim |\O_S(k)|. $$
The fiber of $J$ over generic $C \in |\O_S(k)|$ is thus necessarily
zero-dimensional, and our claim is proved.
It follows that for generic $C \in |kL|$, $\varphi_{|C}$
is of degree 1 onto its image $\varphi(C)$. We can assume $C$ to be
smooth. Then it is the normalization of $\varphi(C)$, and these two curves
have the same geometric genus.

A similar argument shows that for $C \in |\O_S(k)|$ generic, there
cannot exist three pairwise distinct points on $C$ having the same
image under the action of $\varphi$.
So, since $C$ is smooth and $\varphi_{|C}$ is a local imbedding, all
singular points of $\varphi(C)$ occur as the identification of two
distinct points in $C$ by $\varphi$.
We shall now prove that for generic $C \in |\O_S(k)|$, these singular
points are all nodes. 
Write $p : \P(T_S) \rightarrow S$ for the canonical projection of the
projectivized
holomorphic tangent bundle, and consider the incidence variety 
$$ J' \subset |\O_S(k)| \times \P(T_S) \times \P(T_S), $$
defined by the equations
$$ 
(C,u_1,u_2) \in J' \iff
\left\{ \begin{array}{l}
u_1,u_2 \in \P(T_C) \\
p(u_1) \not= p(u_2) \\
\varphi \circ p(u_1) = \varphi \circ p(u_2) \\
\varphi_* u_1 = \varphi_* u_2.
\end{array} \right.$$
It parametrizes the couples of tangent directions of $S$ at two different
points, that are sent by the differential $\varphi_*$ on a couple of
colinear tangent directions at the same point of $S$ (\emph{i.e.} exactly
the situations that yield degenerated singularities on $\varphi(C)$). 
The image of the projection of $J'$ on $\P(T_S) \times \P(T_S)$ is
given by the 
conditions $p(u_1) \not= p(u_2)$, $\varphi \circ p(u_1) = \varphi
\circ p(u_2)$ and $\varphi_* u_1 = \varphi_* u_2$, which yield three
independent equations, so it is of dimension 3. When $k$ is large
enough, the
generic fiber of $J'$ over its projection on $\P(T_S) \times
\P(T_S)$ is of codimension 4 in $|\O_S(k)|$ (see
\ref{condind}). Then   
$$ \dim J' = \dim |\O_S(k)| -1, $$
and the fiber of $J'$ over generic $C \in |\O_S(k)|$ is necessarily
empty, which concludes the proof.

\hfill $\Box$

\begin{rem}\label{condind}
Proposition \ref{action} works as soon as $k \geqslant 4$ when $g
\geqslant 3$, and as soon as $k \geqslant 6$ when $g=2$.
\end{rem}

Indeed, if $g \geqslant 3$ (resp. $g=2$), then
the line bundle $\O_S(k)$ is very ample for $k \geqslant 2$ (resp. $k
\geqslant 3$). This is sufficient  to ensure that two distinct points
of $S$ impose independent conditions on $|\O_S(k)|$ and thus that the
argument concerning $J$ is correct.

Now let $x_1$ and $x_2$ be two distinct points in some projective space
$\P^N$, and $u_1 \in \P(T_{\P^N,x_1})$, $u_2 \in \P(T_{\P^N,x_2})$ be
two tangent directions. As soon as the line defined by $(x_1,u_1)$
(resp. $(x_2,u_2)$) does not pass through $x_2$ (resp. $x_1$), we are
sure that $(x_1,u_1)$ and $(x_2,u_2)$ impose independent conditions on
the linear system of quadrics in $\P^N$. So when $g \geqslant 3$ and
$k \geqslant 4$ (resp. $g=2$ and $k \geqslant 6$), the claim about the
dimension of the generic fiber of $J'$ is true.

\begin{theo}\label{link}
Let $g,l \geqslant 2$ be given. If for $m \in \M^{\circ}_{K3,g}$
generic there exists
a dominant rational map $\varphi_m : S_m \dashrightarrow S_m$
satisfying $\varphi_m^* \O_S(1) \cong \O_S(l)$, 
then for $k$ great enough the universal Severi variety
$\mathcal{V}_{kl,p_a(k)}$ possesses at least two irreducible
components. 
\end{theo}

Theorem~\ref{link} is one of the main results of this article. The
key of the proof is the construction of two irreducible components of
${V}_{kl,p_a(k)}$ for $S$ generic, such that the rational equivalence
class (in $\mathrm{CH}_0(S)$) of the singularity $0$-cycle is constant
for the curves parametrized by the first component, and non constant
for those parametrized by the second component. We give these two
constructions in lemmas \ref{const1} and \ref{const2}.

\begin{lem}\label{const1}
Under the hypotheses of proposition~\ref{action}, and for $k$ large
enough, there exists an irreducible component of ${V}_{kl,p_a(k)}$ on
which the application
$$ C \in {V}_{kl,p_a(k)} \mapsto \mathrm{cl}\left(Z_C\right) \in
\mathrm{CH}_0(S) $$ 
is constant.
\end{lem}
Here ${V}_{kl,p_a(k)}$ is a Severi variety related to the single
surface $S$, and for $C \in {V}_{kl,p_a(k)}$,
$\mathrm{cl}\left(Z_C\right)$ is the rational
equivalence class of the singularity $0$-cycle $Z_C$ of the curve $C$. 

\noindent \textbf{Proof.} 
By proposition~\ref{action}, for $C \in |kL|$ generic and $k$ large
enough, $\varphi(C)$ is an 
irreducible nodal curve in $|klL|$, with geometric genus $p_a(k)$, 
and therefore $\varphi(C) \in V_{kl,p_a(k)}$. 
$V_{kl,p_a(k)}$ is of the expected dimension $p_a(k)$ by
theorem~\ref{dimfamille}, while $|kL|$ is a projective space of
dimension $p_a(k)$. 
So the subset of $V_{kl,p_a(k)}$ parametrizing the images of curves in
$|kL|$ under the action of $\varphi$ is an irreducible component $V'$ of
$V_{kl,p_a(k)}$. 

Let $C$ be a generic curve in $|kL|$, and write $Z_{\varphi(C)}$ for
the $0$-cycle of 
the singular points of its image $\varphi(C)$, seen as a $0$-cycle in
$S$. From the proof of 
proposition~\ref{action}, we know that $\varphi_{|C}: C \rightarrow
\varphi(C)$ is a normalization of $\varphi(C)$. The latter being an
irreducible nodal curve, we have 
$$ 2Z_{\varphi(C)} = K_{\varphi(C)} - (\varphi_{|C})_* K_C,$$
as $0$-cycles in $\varphi_*C$.
This proves that for another generic curve $C' \in |kL|$, the
singularity $0$-cycle $Z_{\varphi(C')}$ of the image $\varphi(C')$ is
rationally equivalent to $Z_{\varphi(C)}$, as $0$-cycles in $S$.  
Indeed, since $C$ and $C'$ are rationally equivalent, the
adjunction formula tells us that $K_C=\left.(K_S+C)\right|_{C}$ and
$K_{C'}=\left.(K_S+C')\right|_{C'}$ are rationally equivalent, as
$0$-cycles on $S$. 
$\varphi(C)$ and $\varphi(C')$ are rationally equivalent
as well, since they both are in $|klL|$, and the adjunction formula
tells us that 
$K_{\varphi(C)}=\left.(K_S+\varphi(C))\right|_{\varphi(C)}$ and
$K_{\varphi(C')}=\left.(K_S+\varphi(C'))\right|_{\varphi(C')}$ are
rationally equivalent, as $0$-cycles on $S$. 

\hfill $\Box$

\begin{rem}
In fact, one gets
$$
\mathrm{cl} \left( Z_{\varphi(C)} \right) = 
\frac{1}{2} k^2 (l^2 -1) L^2 
= \frac{p_a(kl) - p_a(k)}{2g-2} L^2 \in \mathrm{CH}_0(S).
$$
$\delta = p_a(kl) - p_a(k)$ is the number of nodes of $\varphi(C)$ for
$C$ generic.
Using \cite{chowK3}, we get
$$ \mathrm{cl} \left( Z_{\varphi(C)} \right) = 
\delta c_X \in \mathrm{CH}_0(S), $$
where $c_X$ is the rational equivalence class of any point of $S$ that
lie on a rational curve.
\end{rem}

We now construct an irreducible component of $V_{kl,p_a(k)}$, which
parametrizes curves with non constant rational equivalence class for
their singularity $0$-cycles.
\begin{lem}\label{const2}
Under the hypotheses of proposition~\ref{action}, and for $k$ large
enough, there exists an irreducible component of ${V}_{kl,p_a(k)}$ on
which the application
$$ C \in {V}_{kl,p_a(k)} \mapsto \mathrm{cl}\left( Z_C \right) \in
\mathrm{CH}_0(S) $$ 
is non constant.
\end{lem}

\noindent \textbf{Proof.}
By theorem~\ref{dimfamille}, there exists an irreducible family of
dimension $p_a(k)-1$ of irreducible curves $C \in |kL|$, with only one
node as singularity. We write $Z_C$ for the $0$-cycle on $S$ defined
by the singular point of $C$.
$\varphi(C)$ is generically an irreducible curve in 
$|klL|$, with exactly $\delta+1$ nodes as singularities,
$\delta=p_a(kl)-p_a(k)$ (this is proposition~\ref{action}). 
In fact, one has 
$$Z_{\varphi(C)} = \varphi_* Z_C + Z'_{\varphi(C)}$$
as $0$-cycles on $S$, where $Z'_{\varphi(C)}$ is the sum of the
singular points that appear when applying $\varphi$. As $C$
moves, $Z'_{\varphi(C)}$  has
constant rational equivalence class in $\mathrm{CH}_0(S)$, exactly as
in the proof of lemma~\ref{const1}. 

We thus have an irreducible $(p_a(k)-1)$-dimensional family of curves
$\varphi(C) \in |klL|$. For each $\varphi(C)$, we smooth one of the
nodes that are in $Z'_{\varphi(C)}$. 
This eventually gives an irreducible, $p_a(k)$-dimensional family of
irreducible, nodal curves in $|klL|$, with exactly $\delta$ nodes. 
Such a family is an irreducible component $V''$ of
${V}_{kl,p_a(k)}$. 

We claim that the rational equivalence class of the singularity
$0$-cycles $Z_{C'}$ of curves $C'$ parametrized by $V''$ is non
constant. This can be seen by the following simple consideration. For
any points $x,y \in S$, and if $k$ is large enough, we can find a
curve $C \in |kL|$ with a node at $x$ as its only singular point, and
such that $\varphi(C)$ is nodal, with nodes at $\varphi(x)$ and
$y$. Smoothing $y$, we get curves $C'$ in $V''$, with singularity
$0$-cycle 
$$ Z_{C'} = \varphi(x) + Z'_{\varphi(C)} - y. $$
Since $Z'_{\varphi(C)}$ has constant rational equivalence class,
fixing $x$ and letting $y$ move, we see that that the rational
equivalence class $\mathrm{cl} \left( Z_{C'} \right) \in
\mathrm{CH}_0(S)$ cannot be constant.

\hfill
$\Box$

\noindent \textbf{Proof of theorem~\ref{link}}
We write the Stein factorization
$$\xymatrix{
& \mathcal{V}_{k,p} \ar[d] \ar[dl] \\
\widetilde{\mathcal{M}_{K3,g}^{\circ}} \ar[r] &
\mathcal{M}_{K3,g}^{\circ}
}$$
of the projective morphism $\mathcal{V}_{k,p} \rightarrow
\mathcal{M}_{K3,g}^{\circ}$. $\mathcal{V}_{k,p} \rightarrow
\widetilde{\mathcal{M}_{K3,g}^{\circ}}$ is a 
projective morphism with connected fibers, while
$\widetilde{\mathcal{M}_{K3,g}^{\circ}} \rightarrow
\mathcal{M}_{K3,g}^{\circ}$  is finite. A point of
$\widetilde{\mathcal{M}_{K3,g}^{\circ}}$ over $m \in
\mathcal{M}_{K3,g}^{\circ}$ represents a connected component of
$(\mathcal{V}_{k,l})_m$. The monodromy of this morphism thus acts as a
subgroup of the permutation group of the connected components of
fibers of $\mathcal{V}_{k,p} \rightarrow
\mathcal{M}_{K3,g}^{\circ}$. Irreducibility of $\mathcal{V}_{k,p}$
is equivalent to the fact that the monodromy acts transitively on the
components of the 
fibers $V_{k,l}$ (see \emph{e.g.} \cite{severi}). 

If there exists a dominant rational map $\varphi_m : S_m \dashrightarrow
S_m$, satisfying $\varphi_m^* \O_{S_m}(1) \cong \O_{S_m}(l)$ for
generic $m \in \mathcal{M}_{K3,g}^{\circ}$, then we have by lemmas
\ref{const1} and \ref{const2} two irreducible components $V'_m$ and
$V''_m$ of each generic fiber
$(\mathcal{V}_{kl,p_a(k)})_m$ that are algebraically distinguished,
since for curves parametrized by the first one, all singularity
$0$-cycles are rationally equivalent, 
and for curves parametrized by the other one, they are not. It follows
that there exists an open subset 
$$\xymatrix@C=0pt{
&& \widetilde{\mathcal{M}_{K3,g}^{\circ}} \ar[d] \\
U &\subset& \mathcal{M}_{K3,g}^{\circ},
}$$
such that all fibers over $U$ contain at least two points that are
algebraically distinguished. The monodromy cannot
exchange these two points. In particular it does not act transitively,
and $\mathcal{V}_{kl,p_a(k)}$ is not irreducible.

\hfill $\Box$

\section{Properties of a self-rational map on a $K3$
  surface}\label{contraintes}

This section is devoted to the study of a dominant
self-rational map on a given $K3$ surface. The observation of the geometry of
an elimination of indeterminacies gives properties that this map
must satisfy, and which of course restrain
the possibilities for such a self-rational map to exist.
We first get numerical relations between the algebraic and topological
degree that are always valid. We then make further remarks depending
on the complexity of the
elimination of indeterminacies, and give a way to control this
complexity.

The notations are as follows.
$S$ is a generic algebraic $K3$ surface. We assume in
particular that $\Pic(S) = \Z \cdot L$,
where $L$ is effective and satisfies $L^2 = 2g-2$ ($g \in
\mathbf{N}^*$). $\varphi : S \dashrightarrow S$ is a dominant
rational map, and $l$ the positive integer
such that $\varphi^* \O_S(1) \cong \O_S(l)$.
We assume $l>1$. 
We consider an elimination of indeterminacies of
$\varphi$, \emph{i.e.} a commutative diagram 
$$\xymatrix{
\widetilde{S} \ar[d]_{\tau} \ar[dr]^{\widetilde{\varphi}} & \\
S \ar@{-->}[r]_{\varphi} & S,
}$$
where $\tau$ is a finite sequence of blow-ups
$$\widetilde{S} = S_p \overset{\epsilon_p}{\longrightarrow} S_{p-1}
\overset{\epsilon_{p-1}}{\longrightarrow} \cdots
\overset{\epsilon_2}{\longrightarrow} S_1
\overset{\epsilon_1}{\longrightarrow} S_0 = S.$$
We write $F_i$ for the exceptional divisor which appears with
$\epsilon_i$, and $E_i$ for $\epsilon_p^* \circ \cdots \circ
\epsilon_{i+1}^* F_i$ ($1 \leqslant i \leqslant p$).
$(\tau^* L,E_1,\ldots,E_p)$ is an orthogonal basis of $\Pic(\Stilde)$, and
$E_i^2 = -1$ ($1 \leqslant i \leqslant p$).

\subsection{Numerical properties}

We start with a numerical observation coming from Hodge theory.
\begin{prop}\label{carre}
There exists an integer $\lambda$, such that
$$ \deg \varphi = \lambda^2. $$
\end{prop}

\noindent \textbf{Proof.} Let $\omega$ be a global, nowhere vanishing,
holomorphic 2-form on $S$. Since $K_S$ is trivial, and
$K_{\widetilde{S}} = \tau^* K_S + E_1 + \cdots + E_p$, where the
$E_i$'s are exceptional divisors, any global holomorphic 2-form on
$\Stilde$ is a multiple of $\tau^* \omega$. In particular there exists
$\lambda \in \C$ such that $\phitilde^* \omega = \lambda \tau^* \omega$.

We write $\H^2(S,\Q)_{\mathrm{tr.}}$ for the transcendental rational
cohomology of $S$, that is the orthogonal in $\H^2(S,\Q)$ of the
Neron-Severi group $\mathrm{NS}(S)$, with respect to the intersection
form $\left<\ ,\ \right>$. We shall show that $\phitilde^* \eta =
\lambda \tau^* \eta$ for all 
$\eta \in \H^2(S,\Q)_{\mathrm{tr.}}$ 

We clearly have $\H^2(S,\Q)_{\mathrm{tr.}} \cong
\H^2(\Stilde,\Q)_{\mathrm{tr.}}$ \emph{via} $\tau^*$, and since
$\phitilde^*$ sends the transcendental cohomology classes of $S$ to
transcendental cohomology classes in $\Stilde$, there exists a Hodge
structure morphism 
$$ \psi : \H^2(S,\Q)_{\mathrm{tr}} \longrightarrow
\H^2(S,\Q)_{\mathrm{tr}}, $$
such that for all $\eta \in \H^2(S,\Q)_{\mathrm{tr.}}$, one has
$\widetilde{\varphi}^* \eta = \tau^*(\psi (\eta)).$ Now $\omega
\in \H^2(S,\Q)_{\mathrm{tr.}}$, and $\phitilde^* \omega = \lambda
\tau^* \omega$, so the eigenspace $E_{\lambda}$ relative to
$\lambda$ for $\psi$ is non nempty.

Suppose $E_{\lambda} \subset \H^2(S,\C)_{\mathrm{tr}}$ is a
proper subspace for $S$ generic. Since $\lambda \in \C$ is algebraic
over $\Q$, the equations defining $E_{\lambda}$ are contained in the
countable set of equations with coefficients in $\overline{\Q}$. This says
that, when $S$ moves, $\omega \in E_{\lambda}$  
is contained in a countable union of proper linear subspaces of
$\H^2(S,\C)$. 
This contradicts the surjectivity of
the period map for $K3$ surfaces~: its image is an open set of a
projective quadric in $\P\left(\H^2(S,\C)\right)$ (see \emph{e.g.}
\cite{palaiseau}). We
thus have $E_{\lambda} = \H^2(S,\C)_{\mathrm{tr}}$ for $S$ generic. 

$\psi$ acts on $\H^2(S,\Q)_{\mathrm{tr.}}$ as multiplication by
$\lambda$, so $\lambda$ is necessarily a rational number. 
From the two equalities
$$ 
\int_{\Stilde} \phitilde^* \omega \wedge
\phitilde^* \overline{\omega} = \deg(\varphi) \int_S \omega
\wedge \overline{\omega},
$$
and
$$
\int_{\Stilde} \phitilde^* \omega \wedge
\phitilde^* \overline{\omega} = \lambda^2
\int_{\Stilde} \tau^* \omega \wedge \tau^* \overline{\omega} =
\lambda^2 \int_S \omega \wedge \overline{\omega},
$$
we get $\deg(\varphi) = \lambda^2$. Since $\deg \varphi$ is an integer,
and $\lambda$ is a rational number, $\lambda$ is necessarily an
integer.

\hfill $\Box$

We shall now prove a divisibility property involving $\lambda$. To do
so, we need some further notations. Let us write
$$
\phitilde^* L = l \tau^* L - \sum_{1 \leqslant i \leqslant p} \alpha_i
E_i
$$
for some integers $\alpha_1,\ldots,\alpha_p$. Since $\Pic(S) = \Z
\cdot L$, and the $E_i$'s are effective, there also exist
non negative integers $\beta_1,\ldots,\beta_p$, such that for all $i$
$$ \phitilde_* E_i = \beta_i L.$$
Note that by projection formula,
$$
\alpha_i = 
\phitilde_* (\phitilde^* L \cdot E_i) = 
L \cdot \phitilde_* E_i  =
\beta_i L^2 =
(2g-2) \beta_i.
$$

\begin{lem}\label{divise}
$2g-2$ necessarily divides $l-\lambda$.
\end{lem}

\noindent \textbf{Proof.} Let $\lambda$ be as in the proof
of proposition~\ref{carre}. We
have $\phitilde^* \eta' = \lambda \tau^* \eta'$ for all class $\eta'
\in \H^2(S,\Q)_{\mathrm{tr.}}$. We have on the other hand $\phitilde^*
c_1(L) = l \tau^* c_1(L) - (2g-2) \sum_i \beta_i [E_i]$. 

Any $\eta \in \H^2(S,\Q)$ decomposes over $\Q$ into $\eta = \eta'
+ \eta''$, where $\eta' \in \H^2(S,\Q)_{\mathrm{tr.}}$, and $\eta''=
(\left< \eta,c_1(L) \right> / (2g-2) ) c_1(L)$. Then 
\begin{eqnarray*}
\phitilde^* \eta &=& \lambda \tau^* \eta' + l \tau^* \eta'' - \left<
  \eta,c_1(L) \right> \sum_i \beta_i [E_i] \\
&=& \lambda \tau^* \eta + (l-\lambda) \tau^* \eta'' - \left<
  \eta,c_1(L) \right> \sum_i \beta_i [E_i].
\end{eqnarray*}

The intersection product is unimodular, and $c_1(L)$ is indivisible. So
there exists a class $\eta_1 \in \H^2(S,\Z)$, such that 
$$ \left< \eta_1 , c_1(L) \right> = 1. $$
It decomposes over $\Q$ into $\eta_1 = \eta'_1 + \eta''_1$, and the
equality 
$$
(l-\lambda) \tau^* \eta''_1 = \phitilde^* \eta_1 - \lambda \tau^*
\eta_1  + \left< \eta_1,c_1(L) \right> \sum_i \beta_i [E_i]
$$
shows that $$(l-\lambda) \tau^* \eta''_1 = \frac{l-\lambda}{2g-2}
\tau^* c_1(L)$$ is an
integral cohomology class. Since $c_1(L)$ is indivisible, this shows that
$2g-2$ necessarily divides $l-\lambda$.

\hfill $\Box$

We now look more specifically at the geometry of the elimination of
indeterminacies. To have a more accurate description of the situation,
we consider the proper transforms $\hat{F}_i \subset \Stilde$ of the
$F_i$'s, and introduce their intersection tree. Later on, we will call
it the exceptional tree, or the ramification tree (recall that since
$K_S$ is trivial, the ramification divisor of $\phitilde$ and the
total exceptional divisor of $\tau$ are equal, \emph{cf.} proof of 
proposition \ref{action}).
Its vertices are the
$\hat{F}_i$'s, and two vertices are connected if and only if the two
corresponding divisors meet in $\Stilde$. The descendants of a vertex
$\hat{F}_i$ are the vertices situated below $\hat{F}_i$ in the tree,
\emph{i.e.} those corresponding to divisors whose projection by
$\epsilon_{i+1} \circ \cdots \circ \epsilon_p$ is contained in
$F_i$. The depth $m_i$ of a vertex $\hat{F}_i$ is the number of
ancestors of $\hat{F}_i$ in the tree, \emph{i.e.} the number of points
situated above $\hat{F}_i$. The depth of the tree is the maximal depth
of its vertices.

\begin{exem}
The following exceptional tree
$$
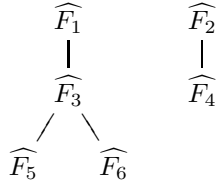
\xymatrix@C=0pt@R=10pt{
          & \hat{F_1} \ar@{-}[d] &           && \hat{F_2} \ar@{-}[d]\\
          & \hat{F_3} \ar@{-}[dr] \ar@{-}[dl] &&           & \hat{F_4} \\
\hat{F_5} &           & \hat{F_6} & \quad &           \\
}$$
is obtained by first blowing up $S$ along two points~; $\hat{F_1}$ and 
$\hat{F_2}$ are the exceptional divisors above these two points. One
          then blows up the resulting surface along one point on
          $\hat{F_1}$, and one point on $\hat{F_2}$. Write $\hat{F_3}$
(resp. $\hat{F_4}$) for the exceptional divisor appearing above the
          blown up point on $\hat{F_1}$ (resp. $\hat{F_2}$). One
          finally blows up along two points of $\hat{F_3}$.
The descendants of $\hat{F_3}$ are $\hat{F_3}$, $\hat{F_5}$ and
          $\hat{F_6}$. Its ancestors are $\hat{F_1}$ and
          $\hat{F_3}$. Its depth is 2. The depth of the tree is 3.
\end{exem}

This being set, $E_i = \epsilon_p^* \circ \cdots \circ \epsilon_{i+1}^*
F_i$ is clearly the sum of all descendants of $\hat{F}_i$ in the tree.
In the above example we have $E_3 = \hat{F_3} + \hat{F_5} +
\hat{F_6}.$
The canonical divisor of $\Stilde$ is $K_{\Stilde} = \epsilon_p^*
\circ \cdots \circ \epsilon_1^* K_S + E_1 + \cdots + E_p$. Since $K_S$
is trivial we have
$$ K_{\Stilde} = \sum_{1 \leqslant i \leqslant p} E_i =  \sum_{1
  \leqslant i \leqslant p} m_i \hat{F}_i. $$ 
It is also the ramification divisor of the map $\phitilde$. 

Let $F$ be an exceptional divisor, such that $\tau$ does not contain
any blow up along a point of $F$ (\emph{i.e.} an exceptional divisor
which appears at the bottom of the exceptional tree). For a suitable
choice of notations, this divisor can be supposed to be $F_p$. 
If $F$ collapses under the action of $\phitilde$, then there
necessarily exists a morphism 
$\phitilde_{p-1} : S_{p-1} \rightarrow S$, and a commutative diagram
$$
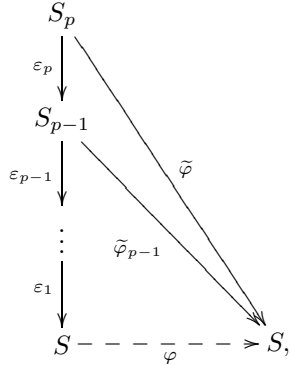
\xymatrix@C=60pt{
S_p \ar[d]_{\epsilon_p} \ar[3,1]^{\phitilde} & \\
S_{p-1} \ar[d]_{\epsilon_{p-1}} \ar[2,1]_{\phitilde_{p-1}} & \\
\vdots \ar[d]_{\epsilon_1} & \\
S \ar@{-->}[r]_{\varphi} & S,
}$$
that is another elimination of indeterminacies of $\varphi$ involving
one less exceptional divisor. 
We may thus assume $\tau$ to be minimal, in
the sense that $\phitilde$ does not contract to a point any exceptional 
divisor which appears at the end of the exceptional tree.

The following equality is obtained simply by computing the
self-intersection $(\phitilde L)^2$. It is the most important relation
between $\deg \varphi$ and $l$. We use the minimality of $\tau$
to show the positivity of the $\beta_i$'s.
\begin{prop}\label{decomp}
The $\beta_i$'s are all positive. 
In addition, the algebraic degree $l$ and the topological degree of
$\varphi$ satisfy   
$$
l^2 = \deg \varphi + (2g-2) \sum_{1 \leqslant i \leqslant p}
\beta_i^2.
$$
\end{prop}

\noindent \textbf{Proof.} 
We have $E_i = \hat{F}_i + \hat{F}_{i_1} + \cdots + \hat{F}_{i_q}$,
  where $\hat{F}_i, \hat{F}_{i_1}, \ldots, \hat{F}_{i_q}$ are all the
  descendants of $\hat{F}_i$ in the exceptional tree. Therefore
$$ \beta_i = \gamma_i + \gamma_{i_1} + \cdots  \gamma_{i_q}, $$
where $\phitilde_* \hat{F}_{i_s} = \gamma_{i_s} L$, $1 \leqslant s
\leqslant q$. The $\gamma_{i_s}$'s are \emph{a priori} non negative
integers. 
$\hat{F}_i$ has at least one descendant $\hat{F}_{i_j}$ at the end of
the exceptional tree. By minimality of $\tau$, $\phitilde$ cannot
contract  $\hat{F}_{i_j}$ to a point, and we have $\gamma_{i_j} \geqslant
1$. Finally
$$ \beta_i \geqslant \gamma_{i_j} >0. $$
We get the relation between $l$ and $\deg \varphi$ simply by computing
in two different ways the self-intersection $(\phitilde^* L)^2$. We
have on the one hand  
$$ (\phitilde^* L)^2 = (\deg \phitilde) L^2 = (\deg \varphi)(2g-2),$$ 
and on the other hand 
$$ (\phitilde^* L)^2 = l^2 (\tau^* L)^2 + \sum_{1 \leqslant i
  \leqslant p} \alpha_i^2 E_i^2 = (2g-2) l^2 - (2g-2)^2 \sum_{1 \leqslant i
  \leqslant p} \beta_i^2,$$
which yields the announced formula.

\hfill $\Box$

We now get the following arithmetic property on the $\beta_i$'s by
some Riemann-Roch computations.

\begin{lem}
$\sum_{1 \leqslant i \leqslant p} \beta_i$ is even.
\end{lem}

\noindent \textbf{Proof.}
We first show that $\phitilde_* \O_{\Stilde}$ is a locally free sheaf
of rank $r := \deg \varphi$. Since it is clearly torsion free, it is
enough to show that any section defined on a punctured open set
$U\setminus\{x_0\}$ extends in a unique way to a section defined over
$U$ (see \cite{reflexive}, lemma 1). So let $U \subset S$ be an open
set, $x_0 \in U$, and $f \in \phitilde_*
\O_{\Stilde}(U\setminus\{x_0\})$. $f$ can be seen as a holomorphic
function on $\Stilde$, defined over
$\phitilde^{-1}(U\setminus\{x_0\})$. If the 
fiber of $\phitilde$ above $x_0$ is a finite set of points, then the
result is clear. Otherwise the fiber contains an irreducible exceptional
curve $F$. $f$ cannot be singular along $F$, since this would give by
restriction a global section of $\O_{\Stilde}(mF)_{|F}$ for some
positive $m$, which is impossible, since $F^2 <0$. So $f$ has only
isolated singularities along $F$, and therefore extends to a
function over $\phitilde^{-1}(U)$.

Now it is an easy consequence of
Grauert's theorem that $R^i \phitilde_* \O_{\Stilde}=0$ for $i >0$.
This gives $\phitilde_! \O_{\Stilde} = \phitilde_*
\O_{\Stilde}$, and  we thus have
\begin{eqnarray*}
ch(\phitilde_! \O_{\Stilde}).td(T_S) &=&
\left( r[S]+ c_1(\phitilde_* \O_{\Stilde}) + \frac{c_1(\phitilde_*
    \O_{\Stilde})^2 -2c_2(\phitilde_* \O_{\Stilde})}{2} \right).([S]+2) 
\\ &=& 
r[S] + c_1(\phitilde_* \O_{\Stilde}) + \left(
\frac{c_1(\phitilde_*\O_{\Stilde})^2 -2c_2(\phitilde_* \O_{\Stilde})}{2}
   + 2r \right).
\end{eqnarray*}
On the other hand, we have
$$
\phitilde_*\left(  ch(\O_{\Stilde}).td(T_{\Stilde}) \right) = 
\phitilde_*\left( [\Stilde] - \frac{1}{2} (E_1 + \cdots + E_p) + 2
\right) = r [S] - \frac{1}{2} \left(\textstyle{\sum_i} \beta_i \right)
L + 2,
$$
so the Grothendieck-Riemann-Roch formula gives
$$
c_1(\phitilde_* \O_{\Stilde}) = - \frac{1}{2} \left(\textstyle{\sum_i}
  \beta_i \right) 
L.$$
Since $L$ is indivisible, the lemma follows.

\hfill $\Box$

\subsection{Complexity of an elimination of indeterminacies}

To motivate the study of the complexity of the elimination of
indeterminacies, we first show that we have further numerical
constraints on $\varphi$ when the elimination of indeterminacies is
not too complicated.
The following numerical property is true under the hypothesis that the 
exceptional tree has depth smaller than $2$.

\begin{prop}\label{amerik}
If the differential $d\phitilde$ does not vanish identically along any
curve of $\Stilde$, then
the topological degree of $\varphi$ satisfies the inequality
$$
\deg \varphi \leqslant 1 + \frac{1}{24}\left[ p + 4(g-1) 
\textstyle{\sum_i}  \beta_i \right].
$$
The condition on the differential is satisfied as soon as the total
depth of the exceptional tree is non greater than 2.
\end{prop}

\noindent \textbf{Proof.} We follow an idea of Amerik, Rovinsky and
Van de Ven (\cite{arv}, see \cite{ba} as well). The fiber bundle
$\Omega^1_S(2)$ is generated by its global sections, so by lemma 1.1 of
\cite{arv} a generic section $\sigma \in \H^0(S,\Omega^1_S(2))$ has
isolated zeroes. With the assumption made on $d\phitilde$ this is also
true for the pull-back section $\phitilde^* \sigma \in
\H^0(\Stilde,\Omega^1_{\Stilde}(2\phitilde^* L)).$ Counting these
zeroes yields the inequality on Chern classes 
$$
c_2\left(\Omega^1_{\Stilde}(2\phitilde^* L)\right) \geqslant (\deg
\varphi) c_2\left(\Omega^1_S(2)\right). 
$$
The left-hand side of this inequality is
$$
c_2(\Omega^1_{\Stilde}) + 
2 \phitilde^* c_1(L) \cdot c_1(\Omega^1_{\Stilde}) 
+ 4 \phitilde^* c_1(L)^2,
$$
and its right-hand side is
$$
\deg \varphi \left[
c_2(\Omega^1_{S}) + 
2 c_1(L) \cdot c_1(\Omega^1_{S}) 
+ 4 c_1(L)^2
\right].
$$
Now $\phitilde^* c_1(L)^2 = (\deg \varphi) c_1(L)^2$, so we get
$$
c_2(\Omega^1_{\Stilde}) + 
2 \phitilde^* c_1(L) \cdot c_1(\Omega^1_{\Stilde}) 
\geqslant 
\deg \varphi \left[
c_2(\Omega^1_{S}) + 
2 c_1(L) \cdot c_1(\Omega^1_{S}) \right],
$$
that is
$$
\chi_{\mathrm{top}}(\Stilde) + 2 \phitilde^* L \cdot K_{\Stilde} 
\geqslant 
\deg \varphi \left[ \chi_{\mathrm{top}}(S) + 2 L \cdot K_S  \right],
$$
where $\chi_{\mathrm{top}}$ denotes the topological Euler-Poincar{\'e}
characteristic, that is the alternated sum of the Betti numbers. It is
24 for all $K3$ surfaces. $\Stilde$ is obtained from $S$ by successively
blowing up along $p$ points so $\chi_{\mathrm{top}}(\Stilde)=24+p$. We
also have $K_S=0$, and
$$
\phitilde^* L \cdot K_{\Stilde} =  (l \tau^* L - 
\textstyle{\sum_i} \alpha_i E_i) \cdot (E_1+\cdots+E_p) =
\textstyle{\sum_i} \alpha_i. 
$$
We eventually get
$$
24+p + 2 \textstyle{\sum_i} \alpha_i \geqslant
24 \deg \varphi,
$$
which yields the desired inequality with the relations $\alpha_i =
(2g-2)\beta_i$.

Now suppose $d\phitilde$ vanishes identically along a curve $C$ in
$\Stilde$. Then $C$ necessarily collapses under the action of
$\phitilde$, and it appears with multiplicity at least 2 in its
ramification divisor. Indeed, let $f$ be some local
equation for $C$. If 
$d\phitilde$ vanishes with order $\mu$ along $C$, then it writes
$$d\phitilde = \left( \begin{array}{cc}
f^{\mu} g_{11} & f^{\mu} g_{12} \\
f^{\mu} g_{21} & f^{\mu} g_{22} 
\end{array} \right)$$
in some local holomorphic coordinate system, with the $g_{ij}$
holomorphic, 
and $\bigwedge^2 d\phitilde$ vanishes with order $2\mu$ along $C$. If
the total depth of the exceptional tree is less than 2, the only
curves which appear with multiplicity greater than 2 in the
ramification divisor are at the end of the tree, and cannot be
contracted to a point by $\phitilde$ by minimality of the elimination of
indeterminacies. So in this case, $d\phitilde$ does not vanish
identically along any curve of $\Stilde.$

\hfill $\Box$

The first step towards a control of the complexity of the elimination
of indeterminacies
is made with the following basic remark. It shows that the depth of
the exceptional tree is controlled by the topological
degree.

\begin{prop}\label{depth}
(i) The depth $m$ of the exceptional tree always satisfy 
$$ m \leqslant \deg \varphi - 2.$$ 
(ii) If the tree has two connected components of depths $m'$ and
$m''$, then 
$$m'+m'' \leqslant \deg \varphi -2.$$
In particular, if one has equality in (i), then the tree only has one
connected component.
\end{prop}

\noindent \textbf{Proof.}
\emph{(i)} 
Since the ramification divisor of $\phitilde$ is $\sum_i m_i
\hat{F}_i$, it is clear that 
$$ m = \max m_i \leqslant \deg \varphi -1.$$
Now suppose there exists an irreducible exceptional curve $F$ that has
depth $\deg \varphi - 1$ in the exceptional tree. Then it is at the
end of the tree, and therefore is not contracted. $F$ appears in the
ramification divisor with multiplicity $\deg \varphi - 1$, and thus
$$ \phitilde^{-1} \left( \phitilde(F) \right) = 
(\deg \varphi) F + E,$$
where $E$ is contracted by $\phitilde$. In particular, $E$ is
exceptional for $\tau$ as $F$ is. It follows that
$\phitilde^{-1} \left( \phitilde(F) \right)$ is supported on the
exceptional divisor of $\tau$, which implies that it has negative
self-intersection. This contradicts the fact that
$$ \phitilde^{-1} \left( \phitilde(F) \right)^2 = 
(\deg \varphi) \phitilde(F)^2 > 0. $$

\emph{(ii)}
If the tree has two connected components of depths $m'$ and
$m''$, then we have two irreducible exceptional curves $F'$ and $F''$
of depths $m'$ and $m''$, that are not contracted, and that do not
meet in $\Stilde$. The image
curves $\phitilde(F')$ and $\phitilde(F'')$ intersect in $S$, because
their images have their class proportional to $c_1(L)$. Let $x$
be an intersection point. There are at least two distinct points $x'
\in F'$ 
and $x'' \in F''$ in $\phitilde^{-1}(x)$. Since $F'$ and $F''$ appear
with multiplicities $m'$ and $m''$ in the ramification divisor of
$\phitilde$, $x'$ and $x''$ appear with multiplicities $m'+1$ and
$m''+1$ in $\phitilde^{-1}(x)$. This implies
$$ m' + m'' + 2 \leqslant \deg \phitilde. $$

\hfill $\Box$

\begin{rem}
In fact, (ii) can be extended as follows~: if there exist two
distinct curves
$F'$ and $F''$ at the end of the exceptional tree, which have depths
$m'$ and $m''$, then 
$$ m' + m'' + 2 \leqslant \deg \phitilde. $$
\end{rem}

Now the following result gives control on another aspect of the
elimination of indeterminacies, namely the number of blown-up points
on $S$. 
It says that in case all irreducible exceptional curves are
disjoint, then their number is bounded from above. The
hypothesis is equivalent to the fact that we can eliminate
the indeterminacies of $\varphi$ by the single blow up of finitely
many distinct points on $S$.

\begin{prop}\label{width}
If the exceptional tree has depth 1, then 
$$
p \leqslant 8(\deg \varphi -1).
$$
\end{prop}

\vspace{0.2cm}\noindent \textbf{Proof.}
In this case, the ramification divisor of $\phitilde$ is $E_1+\cdots
+E_p$, where the $E_i$'s are disjoint $\P^1$'s, and by minimality of
the elimination of 
indeterminacies, $\phitilde$ does not contract any
of them. So the differential $d\phitilde$ is surjective, and we have
an exact sequence 
$$
0 \rightarrow 
\phitilde^* \Omega_S \rightarrow
\Omega_{\Stilde} \rightarrow 
\textstyle{\bigoplus_i} L_i \rightarrow 0,
$$
where each $L_i$ is a line bundle on the exceptional curve $E_i$. This
gives
\begin{eqnarray*}
c_2\left( \Omega_{\Stilde} \right) &=& 
c_2 \left( \phitilde^* \Omega_S \right) + c_1 \left( \phitilde^*
  \Omega_S \right) \cdot c_1 \left( \textstyle{\bigoplus_i} L_i
\right) + c_2 \left( \textstyle{\bigoplus_i} L_i \right)\\
&=&
(\deg \varphi)c_2\left( \Omega_S \right) + \textstyle{\sum_i}
c_2(L_i). 
\end{eqnarray*}
By restriction, we get on each $E_i$ an exact sequence
$$
0 \rightarrow 
K_i \rightarrow 
{\Omega_{\Stilde}}_{|E_i} \rightarrow 
L_i \rightarrow 0,
$$
where $K_i$ is a line bundle on $E_i$. 
We have a map $K_i \rightarrow
\Omega_{E_i}$, given by the composition
$$\xymatrix@=15pt{
& {\Omega_{\Stilde}}_{|E_i} \ar[d] \\
K_i \ar[r] \ar[ur]^{(d\phitilde)^t} & \Omega_{E_i}.
}$$
Since $\phitilde$ is ramified along $E_i$, the map $K_i \rightarrow 
I_i/I_i^2$ induced by $(d\phitilde)^t$ is zero (here $I_i \subset
\O_{\Stilde}$ is the ideal sheaf of $E_i$). This shows that the
above map $K_i \rightarrow \Omega_{E_i}$ is an injection, and thus
that
$$
\deg K_i \leqslant \deg \Omega_{E_i} = -2
$$
(as line bundles on $E_i \cong \P^1$).
On the other hand, $\deg({\Omega_{\Stilde}}_{|E_i}) = -1$ by the
conormal exact sequence, so one has
$$
\deg(L_i)= \deg({\Omega_{\Stilde}}_{|E_i}) -
\deg(K_i) \geqslant 1.
$$
We write $d_i$ for $\deg(L_i)$ (\emph{i.e.} $L_i=\O_{E_i}(d_i)$). The
restriction exact sequence 
$$
0 \rightarrow
\O_{\Stilde}\left(-(d_i+1)E_i\right) \rightarrow
\O_{\Stilde}(-d_iE_i) \rightarrow
\O_{E_i}(d_i) \rightarrow 0
$$
gives the two relations
$$\left\{ \begin{array}{rcl}
c_1(L_i) -(d_i+1)E_i &=& -d_iE_i \\
c_2(L_i) -(d_i+1) E_i \cdot c_1(L_i) &=& 0, 
\end{array} \right. $$
and therefore
$$
c_2(L_i) = -d_i - 1  \leqslant -2.
$$
So eventually 
\begin{eqnarray*}
c_2\left( \Omega_{\Stilde} \right) &=&
(\deg \varphi)c_2\left( \Omega_S \right) + \textstyle{\sum_i} c_2(L_i)
\\ & \leqslant & 24 \deg \varphi -2p,
\end{eqnarray*}
and since one knows that $c_2\left( \Omega_{\Stilde} \right) = 24 +p$,
one gets the announced inequality.

\hfill $\Box$

In case $\deg \varphi = 4$, propositions \ref{depth} and \ref{width}
work very well. The following proposition sums up what they learn us
in this case.

\begin{exem}\label{deg4}
If $\deg \varphi = 4$, then either the exceptional tree has depth $1$
and there are less than $24$ blown-up points, or it has depth $2$ and
there is only one blown up point.
\end{exem}

Note that example~\ref{deg4} shows that if $\deg \varphi = 4$, then
the hypotheses of proposition~\ref{amerik} are always satisfied.
To conclude, we compute the first possible couples
$(\deg \varphi,l)$ for which there could actually be a self-rational
map $\varphi$, according to all numerical properties gathered
above. Recall that $K3$ surfaces of genera $2$, $3$, $4$ and $5$ are
respectively double covers of $\P^2$, quartics in $\P^3$, complete
intersections of a cubic and a quadric in $\P^4$, and complete
intersections of three quadrics in $\P^5$.

\begin{exem}
For $\deg \varphi = 4$, the first possible values of $l$ possible are
given by \\ 
\begin{center}\begin{tabular}{|c|c|c|c|c|}
\hline $g$ & $2$ & $3$ & $4$ & $5$ \\
\hline $l$ & $6,8,10,\ldots$ & $6,10,14,\ldots$ & $8,10,14,\ldots$ &
$6,10,14,\ldots$ \\ 
\hline
\end{tabular} .\end{center}
For $\deg \varphi=9$, we get
\begin{center}\begin{tabular}{|c|c|c|c|c|}
\hline $g$ & $2$ & $3$ & $4$ & $5$ \\
\hline $l$ & $5,7,9,\ldots$ & $5,7,9,\ldots$ & $9,15,21,\ldots$ &
$5,11,13,19,\ldots$ \\ 
\hline
\end{tabular} .\end{center}
\end{exem}


\begin{thebibliography}{BHPV04}

\bibitem[AmCa05]{ac}
Ekaterina Amerik and Fr{\'e}d{\'e}ric Campana.
\newblock Fibrations m{\'e}romorphes sur certaines vari{\'e}t{\'e}s {\`a} fibr{\'e}
  canonique trivial.
\newblock Prepublication math.AG/0510299, 2005.

\bibitem[ARV99]{arv}
E.~Amerik, M.~Rovinsky, and A.~Van~de Ven.
\newblock A boundedness theorem for morphisms between threefolds.
\newblock {\em Ann. Inst. Fourier}, 49:405--415, 1999.

\bibitem[ArCo81]{footnote}
Enrico Arbarello and Maurizio Cornalba.
\newblock Footnotes to a paper of {B}eniamino {S}egre.
\newblock {\em Math. Ann.}, 256(3):341--362, 1981.

\bibitem[Bar77]{reflexive}
Wolf Barth.
\newblock Some properties of stable rank-{$2$} vector bundles on {${\bf
  P}\sb{n}$}.
\newblock {\em Math. Ann.}, 226(2):125--150, 1977.

\bibitem[BHPV04]{barth}
W.~Barth, K.~Hulek, C.~Peters, and A.~Van~de Ven.
\newblock {\em Compact complex surfaces}.
\newblock Ergebnisse der Mathematik und ihrer Grenzgebiete. Springer-Verlag,
  2004.

\bibitem[Bea78]{beauville}
Arnaud Beauville.
\newblock {\em Surfaces alg{\'e}briques complexes}.
\newblock Number~54 in Ast{\'e}rique. Soci{\'e}t{\'e} Math{\'e}matique de France, 1978.
\newblock English traduction available at Cambridge University Press.

\bibitem[Bea01]{ba}
Arnaud Beauville.
\newblock Endomorphisms of hypersurfaces and other manifolds.
\newblock {\em Math. Res. Notices}, (1):53--58, 2001.

\bibitem[BD85]{BeauDo}
Arnaud Beauville and Ron Donagi.
\newblock La vari{\'e}t{\'e} des droites d'une hypersurface cubique de dimension 4.
\newblock {\em C. R. Acad. sci. Paris, S{\'e}rie 1}, 301(14):703--706, 1985.

\bibitem[BV04]{chowK3}
Arnaud Beauville and Claire Voisin.
\newblock On the {C}how ring of a {$K3$} surface.
\newblock {\em J. Algebraic Geom.}, 13(3):417--426, 2004.

\bibitem[BT00]{bt}
Fedor~A. Bogomolov and Yuri Tschinkel.
\newblock Density of rational points on elliptic {$K3$} surfaces.
\newblock {\em Asian J. Math.}, 4(2):351--368, 2000.

\bibitem[Can05]{cantat}
Serge Cantat.
\newblock Syst{\`e}mes dynamiques polynomiaux.
\newblock M{\'e}moire d'habilitation {\`a} diriger des recherches, 2005.

\bibitem[Che99]{chen}
Xi~Chen.
\newblock Rational curves on {$K3$} surfaces.
\newblock {\em J. Alg. Geom.}, 8(2):245--278, 1999.

\bibitem[CC99]{ciliberto}
Luca Chiantini and Ciro Ciliberto.
\newblock On the {S}everi varieties of surfaces in $\mathbf{P}^3$.
\newblock {\em J. Alg. Geom.}, 8(1):67--83, 1999.

\bibitem[Fuj02]{fujimoto}
Yoshio Fujimoto.
\newblock Endomorphisms of smooth projective 3-folds with non-negative
  {K}odaira dimension.
\newblock {\em Publ. Res. Inst. Math. Sci.}, 38(1):33--92, 2002.

\bibitem[FN05]{fuji-naka}
Yoshio Fujimoto and Noboru Nakayama.
\newblock Compact complex surfaces admitting non-trivial surjective
  endomorphisms.
\newblock {\em Tohoku Math. J. (2)}, 57(3):395--426, 2005.

\bibitem[GLS00]{gls}
G.-M. Greuel, C.~Lossen, and E.~Shustin.
\newblock Castelnuovo function, zero-dimensional schemes and singular plane
  curves.
\newblock {\em J. Alg. Geom.}, 9(4):663--710, 2000.

\bibitem[Har86]{severi}
Joe Harris.
\newblock On the {S}everi problem.
\newblock {\em Invent. Math.}, 84(3):445--461, 1986.

\bibitem[HT06]{ht}
Brendan Hassett and Yuri Tschinkel.
\newblock Potential density of rational points for {$K3$} surfaces over
  function fields.
\newblock Prepublication math.AG/0604222, to appear in \emph{Amer. J. Math.},
  2006.

\bibitem[Kei03]{keilen}
Thomas Keilen.
\newblock Irreducibility of equisingular families of curves.
\newblock {\em Trans. Amer. Math. Soc.}, 355(9):3485--3512, 2003.

\bibitem[Kol96]{kollar}
Janos Koll{\'a}r.
\newblock {\em Rational curves on algebraic varieties}.
\newblock Ergebnisse der Mathematik und ihrer Grenzgebiete. Springer-Verlag,
  1996.

\bibitem[Nak02]{nakayama}
Noboru Nakayama.
\newblock Ruled surfaces with non-trivial surjective endomorphisms.
\newblock {\em Kyushu J. Math.}, 56(2):433--446, 2002.

\bibitem[Pal85]{palaiseau}
S{\'e}minaire Palaiseau.
\newblock {\em G{\'e}om{\'e}trie des surfaces {$K3$}~: modules et p{\'e}riodes}.
\newblock Number 126 in Ast{\'e}risque. Soci{\'e}t{\'e} math{\'e}matique de France,
  1985.

\bibitem[SD74]{saintdo}
Bernard Saint-Donat.
\newblock Projective models of {$K3$} surfaces.
\newblock {\em Amer. J. Math.}, 96:602--639, 1974.

\bibitem[Tan82]{tannenbaum}
Allen Tannenbaum.
\newblock Families of curves with nodes on {$K3$} surfaces.
\newblock {\em Math. Ann.}, 260(2):239--253, 1982.

\bibitem[Voi02]{voisin}
Claire Voisin.
\newblock {\em Th{\'e}orie de Hodge et g{\'e}om{\'e}trie alg{\'e}brique complexe}.
\newblock Number~10 in Cours sp{\'e}cialis{\'e}s. Soci{\'e}t{\'e} Math{\'e}matique de
  France, 2002.
\newblock English traduction available at Cambridge University Press.

\bibitem[Voi03]{harvard}
Claire Voisin.
\newblock On some problems of {K}obayashi and {L}ang~; algebraic approaches.
\newblock {\em Current developments in Mathematics}, (1):53--125, 2003.

\bibitem[Voi04]{Kcorresp}
Claire Voisin.
\newblock Intrinsic pseudo-volume forms and {$K$}-correspondences.
\newblock In {\em The Fano conference}, pages 761--792. Univ. Torino, Turin,
  2004.

\bibitem[Zar82]{zariski}
Oscar Zariski.
\newblock Algebraic systems of plane curves.
\newblock {\em Am. J. Math.}, 104(1):209--226, 1982.

\end{thebibliography}
\end{document}